\documentclass[11pt]{article}
\usepackage{amsfonts}
\usepackage{bbm}
\usepackage{mathrsfs}
\usepackage{amssymb}
\usepackage{url}

\title{Generalized Down-up Algebras Revisited\\ from A Viewpoint of Gr\"obner Basis Theory}

\author{Rabigul Tuniyaz\thanks{Project supported by
the National Natural Science Foundation of China (11861061).}\\
{\small Department of Mathematics, School of  Science}\\
{\small Xinjiang Institute of Science and Technology}\\
{\small  Akesu, 843100, Xinjiang, China}\\
{\small~~~}\\
{\small Gulshadam Yunus\thanks{Project supported by
the National Natural Science Foundation of China (12061068).}}\\
{\small College of Mathematics and System Sciences}\\
{\small Xinjiang University}\\
{\small Urumqi, 830046, Xinjiang, China}}

\date{}

\begin{document}
\maketitle
\begin{center}
\begin{minipage}{120mm}
{\small {\bf Abstract.} The so called generalized down-up algebras are revisited from a viewpoint of Gr\"obner basis theory. Particularly it is  shown explicitly that generalized down-up algebras are solvable polynomial algebras (provided $\lambda\omega\ne 0$), and by means of homogeneous Gr\"obner defining relations, the associated graded structures of generalized down-up algebras, namely the associated graded algebras, Rees algebras, and the homogenized algebras of generalized down-up algebras, are explored comprehensively.}
\end{minipage}\end{center} {\parindent=0pt\par

}
\vskip -.5truecm

\renewcommand{\thefootnote}{\fnsymbol{footnote}}
\setcounter{footnote}{-1}


\def\NZ{\mathbb{N}}
\def\QED{\hfill{$\Box$}}
\def \r{\rightarrow}
\def\mapright#1#2{\smash{\mathop{\longrightarrow}\limits^{#1}_{#2}}}

\def\v5{\vskip .5truecm}
\def\OV#1{\overline {#1}}
\def\hang{\hangindent\parindent}
\def\textindent#1{\indent\llap{#1\enspace}\ignorespaces}
\def\item{\par\hang\textindent}

\def\LH{{\bf LH}}\def\LM{{\bf LM}}\def\LT{{\bf
LT}}\def\KX{K\langle X\rangle} \def\KS{K\langle X\rangle}
\def\B{{\cal B}} \def\LC{{\bf LC}} \def\G{{\cal G}} \def\FRAC#1#2{\displaystyle{\frac{#1}{#2}}}
\def\SUM^#1_#2{\displaystyle{\sum^{#1}_{#2}}}
\def\T{\widetilde}

\section*{1. Introduction}

Let $K$ be a field, and $U(\textsf{sl}(2,K))$ the enveloping algebra of the 3-dimensional Lie algebra $\textsf{sl}(2,K)=Kx\oplus Ky\oplus Kz$ defined by the bracket product: $[x,y]=z$, $[z,x]=2x$, $[z,y]=-2y$. Due to their importance in physics, particularly in quantum group theory, deformations of $U(\textsf{sl}(2,K))$ over $K$ have been studied in great generality in various contexts (e.g., [24, 12, 13]), among which, the down-up algebras introduced in [1, 2] and the generalized down-up algebras introduced in [5] have covered a large class of deformations of $U(\textsf{sl}(2,K))$, and many basic structural properties of down-up algebras and generalized down-up algebras have been explicitly established respectively in pure ring theoretic methods, methods of representation theory, and geometric methods (e.g., [4, 3, 8, 11, 9). However, it seems that there has been no an unified approach to the investigation of generalized down-up algebras by using the Gr\"obner basis methods, though in [16, 17]  the Gr\"obner defining relations of more general algebras similar to  $U(\textsf{sl}(2,K))$ (including generalized down-up algebras) have been established.  In this note, we revisit generalized down-up algebras via their  Gr\"obner defining relations, recapture most of the known structural properties of such algebras, show explicitly that if $\lambda\omega\ne0$ then generalized down-up algebras are all solvable polynomial algebras  in the sense of [10] (thereby they are equipped with an effective Gr\"obner basis theory),  and by means of the homogeneous Gr\"obner defining relations, the associated graded structures of generalized down-up algebras, namely the associated graded algebras, Rees algebras, and the homogenized algebras of generalized down-up algebras, are explored comprehensively.\v5

Throughout this note, $K$ denotes an algebraically closed field, $K^*=K-\{0\}$, and all $K$-algebras considered are associative with multiplicative identity 1. If $S$ is a nonempty subset of an algebra $A$, then we write $\langle S\rangle$ for the two-sided ideal of $A$ generated by $S$.\par

Moreover, for the reader's convenience of understanding the Auslander regularity and the Cohen-Macaulay property of an algebra in Section 4 and Section 5, let us recall that a finitely generated algebra $A$ of finite global dimension $n$ is said to \par
(a) be {\it Auslander regular} if for every finitely generated left $A$-module $M$, every integer $j\ge 0$ and every (right) $A$-submodule $N$ of Ext$^j_A(M,A)$ we have that $j(N)\ge j$, where $j(N)$ is the grade number of $N$ which is the least integer $i$ such that Ext$^i_A(M,A)\ne 0$;\par

(b) satisfy the {\it Cohen-Macaulay property} if for every finitely generated left $A$-module $M$ we have the equality: GK.dim$M+j(M)=n$, where GK.dim denotes the Gelfand-Kirollov dimension. \par

Concerning the Auslander regularity and the Cohen-Macaulay property of filtered rings and graded rings, in particular, of the generalized down-up algebras and their associated graded structures, one is referred to [13, 18, 20].\v5

\section*{2. The Gr\"obner defining relations of Generalized down-up algebras}
In this section we review, in a little more detail,  the introduction of generalized down-up algebras in the sense of [5], so as to strengthen the connection of such algebras with other important algebras. Moreover, by referring to [16, 17, 5], we conclude that the set of defining relations of a generalized down-up algebra forms a  Gr\"obner basis in the sense of [23].\v5

Let $\KS =K\langle X_1,X_2,X_3\rangle$ be the free
$K$-algebra generated by $X=\{ X_1,X_2,X_3\}$. Consider the algebra $A=\KS
/\langle\G\rangle$ with $\G$ consisting of
$$(*)~\left\{\begin{array}{l} g_{31}=X_3X_1-\lambda X_1X_3+\gamma X_3,\\
g_{12}=X_1X_2-\lambda X_2X_1+\gamma X_2,\\
g_{32}=X_3X_2-\omega X_2X_3+f(X_1),\end{array}\right.$$ where $\lambda
,\gamma ,\omega\in K$, and $f(X_1)$ is a polynomial in the variable
$X_1$. In [5], this algebra is referred to as a {\it generalized down-up algebra}, and is denoted by $L(f, \lambda,\omega ,\gamma)$.\v5

{\bf Convention} For the purpose of this note and for saving notation, throughout this note we will use $A$, with the presentation $A=\KS /\langle\G\rangle$, to denote a generalized down-up algebra with the set of defining relations $\G$ as described above. \v5

As one may see from the literature, or as illustrated below, the algebras $A$ defined above are mainly stemming from two topics:\par

(1) The investigation of deformations of the enveloping algebra $U(\textsf{sl}(2,K))$ of the 3-dimensional Lie algebra $\textsf{sl}(2,K)$. \vskip 6pt

(a)  With $\lambda
=\omega =1$, $\gamma =2$ and $f(X_1)=-X_1$ in $(*)$, it is clear that $A=U(\textsf{sl}(2,K))$.\vskip 6pt

(b) In [24], Smith introduced a family of algebras $A$ similar to $U(\textsf{sl}(2,K))$, that is, $\G$ consists of
$$\begin{array}{l} g_{31}=X_1X_3-X_3X_1=X_3,\\ g_{12}=X_1X_2-X_2X_1=-X_2,\\ g_{32}=X_3X_2-X_2X_3=f(X_1).\end{array}$$\par

(c)  Let $\zeta\in K^*$. With $\lambda =\zeta^4$, $\omega =\zeta^2$,
$\gamma =-(1+\zeta^2)$, $a=0=c$, and $b=-\zeta$, the algebra $A$ coincides with Woronowicz's deformation of $U(\textsf{sl}(2,K))$ which was introduced in the noncommutative differential calculus [26].\vskip 6pt

(d) If $f(X_1)=bX_1^2+X_1$ and $\lambda\gamma \omega b\ne 0$, then $A$ coincides with Le Bruyn's conformal $\textsf{sl}_2$ enveloping algebra [13, Lemma 2] which provides a special family of Witten's deformation
of $U(\textsf{sl}(2,K))$ in quantum group theory [28].\vskip 6pt

(2) The investigation of down-up algebras in the sense of [1, 2]. \par
(e) Let $K\langle X_1, X_2\rangle$ be the free $K$-algebra generated by $\{ X_1, X_2\}$.  The down-up algebra $A(\alpha ,\beta ,\gamma )$, in the sense of [1, 2], was introduced in  the study of algebras generated by the down and up operators on a differential or uniform partially ordered set (poset),  that is $A(\alpha ,\beta ,\gamma )=K\langle X_1, X_2\rangle /\langle S\rangle$ with $S$ consisting of
$$\begin{array}{l} f_1=X_1^2X_2-\alpha X_1X_2X_1-\beta X_2X_1^2-\gamma X_1,\\
f_2=X_1X_2^2-\alpha X_2X_1X_2-\beta X_2^2X_1-\gamma X_2,\end{array}$$
where $\alpha$, $\beta$, $\gamma\in K$. If, in the foregoing definition of a generalized down-up algebra $A$, the polynomial $f(X_1)$ has degree one, then all down-up algebras $A(\alpha ,\beta ,\gamma )$ are retrieved for suitable choices of the parameters of $A$ [9, 4, 5, 3], that is, each down-up algebra is isomorphic to some generalized down-up algebra $A$. For more detailed argumentation on this result, see [3, Lemma 1.1]. \v5

At this stage, also let us point out that the down-up algebras have been connected with many more important algebras, for instance, the generalized Weyl algebra (see [5]), the hyperbolic rings (see [11, Proposition 3.0.1]), and  the parafermionic (parabosonic) algebra which is closely related to the cubic Artin-Schelter regular algebras (see [7]), namely the natural action of $GL(2)$ on the parafermionic algebra for $D=2$ extends as an action of the quantum group $GL_{p,q}(2)$ on the generic cubic Artin-Schelter regular algebra of type $S_1$ with the defining relations
$$\begin{array}{l} g_1=X_2X_1^2+qrX_1^2X_2-(q+r)X_1X_2X_1,\\
g_2=X_2^2X_1+qrX_1X_2^2-(q+r)X_2X_1X_2,\end{array}$$
where $q,r\in\mathbb{C}$.\v5

Now, let $A$ be a generalized down-up algebra. Then, in either of the following two cases:{\parindent=1.25truecm\par
\item{(a)} deg$f(X_1)=n\le 2$, $X_1$, $X_2$, and $X_3$ are all assigned  the degree $1$;\par
\item{(b)} deg$f(X_1)=n\ge 1$,  $X_1$ is assigned the degree 1, but $X_2$ and $X_3$ are all assigned  the degree $n$,\par}{\parindent=0pt
the set $\G =\{ g_{31},g_{12},g_{32}\}$ of defining relations of $A$ forms a Gr\"obner basis for the ideal $\langle\G\rangle$ in the sense of [23], where in both cases the monomial
ordering used on $\KS$ is the graded lexicographic ordering
$$X_2\prec_{grlex} X_1\prec_{grlex}X_3.$$ Thereby $A$ has the PBW $K$-basis $\B =\{ a_2^ia_1^ja_3^{\ell}~|~(i,j,\ell)\in\NZ^3\}$, where $a_k$ is the coset represented by $X_k$ in $A$, $1\le k\le 3$. One may refer to ([17, Ch.4, Section 3], [5, Theorem 2.1]) for detailed argumentations, though the notion of a Gr\"obner basis is not obviously used by [5]. For convenience of later usage, we especially record this result  here.}\v5

{\bf Proposition 2.1} With notation as above, the following statements hold.\par
(i) In both the cases (a) and (b) above, $\G$ is a Gr\"obner basis of the ideal $\langle\G\rangle$ with respect to the monomial ordering $X_2\prec_{grlex}X_1\prec_{grlex}X_3$ on $\KS$.\par
(ii) The generalized algebra $A=\KS /\langle\G\rangle$ has the PBW $K$-basis $$\B =\{ a_2^ia_1^ja_3^{\ell}~|~(i,j,\ell)\in\NZ^3\},$$ where $a_k$ is the coset represented by $X_k$ in $A$, $1\le k\le 3$.\par\QED\v5

{\bf Remark} The reason that in the above argumentation we distinguish the cases (a) and (b), is to obtain Theorem 4.2(i), Theorem 4.3(ii) and Theorem 5.3(iv) in later Section 4 and Section 5 respectively.\v5

\section*{3. Generalized down-up algebras are solvable polynomial algebras provided $\lambda\omega\ne 0$}
In an example of [17, P.154], without proof in detail it is concluded that if $\lambda\omega\ne 0$, then the corresponding generalized down-up algebras are all solvable polynomial algebras in the sense of [10]. It means that the Gr\"obner basis theory and computational methods for solvable polynomial algebras [10, 16, 17] may completely applied to generalized down-up algebras. Instead of going to make a long story about this topic, in this section we just give a detailed argumentation of the conclusion mentioned above, and from this fact we re-derive two basic structural properties of generalized down-up algebras given in [5].\v5

First recall from [10, 21, 16, 17] the following definitions. Suppose that a finitely generated  $K$-algebra $A=K[a_1,\ldots ,a_n]$ has the PBW $K$-basis $\B =\{ a^{\alpha}=a_{1}^{\alpha_1}\cdots
a_{n}^{\alpha_n}~|~\alpha =(\alpha_1,\ldots ,\alpha_n)\in\NZ^n\}$, and that $\prec$ is a total ordering on $\B$. Then every nonzero element $f\in A$ has a unique expression
$$\begin{array}{rcl} f&=&\lambda_1a^{\alpha (1)}+\lambda_2a^{\alpha (2)}+\cdots +\lambda_ma^{\alpha (m)},\\
&{~}&\hbox{such that}~a^{\alpha (1)}\prec a^{\alpha
(2)}\prec\cdots \prec a^{\alpha (m)},\\
&{~}&\hbox{where}~ \lambda_j\in K^*,~a^{\alpha
(j)}=a_1^{\alpha_{1j}}a_2^{\alpha_{2j}}\cdots a_n^{\alpha_{nj}}\in\B
,~1\le j\le m.
\end{array}$$
Noticing that  elements of $\B$ are conventionally called {\it monomials}\index{monomial}, the {\it leading monomial of $f$} is defined as $\LM
(f)=a^{\alpha (m)}$, the {\it leading coefficient of $f$} is defined
as $\LC (f)=\lambda_m$, and the {\it leading term of $f$} is defined
as $\LT (f)=\lambda_ma^{\alpha (m)}$.\v5

{\bf Definition 3.1}  Suppose that the $K$-algebra
$A=K[a_1,\ldots ,a_n]$ has the PBW basis $\B$. If $\prec$ is a
total ordering on $\B$ that satisfies the following three
conditions:{\parindent=1.35truecm\par

\item{(1)} $\prec$ is a well-ordering (i.e., every nonempty subset of $\B$ has a minimal element);\par

\item{(2)} For $a^{\gamma},a^{\alpha},a^{\beta},a^{\eta}\in\B$, if $a^{\gamma}\ne 1$, $a^{\beta}\ne
a^{\gamma}$, and $a^{\gamma}=\LM (a^{\alpha}a^{\beta}a^{\eta})$,
then $a^{\beta}\prec a^{\gamma}$ (thereby $1\prec a^{\gamma}$ for
all $a^{\gamma}\ne 1$);\par

\item{(3)} For $a^{\gamma},a^{\alpha},a^{\beta}, a^{\eta}\in\B$, if
$a^{\alpha}\prec a^{\beta}$, $\LM (a^{\gamma}a^{\alpha}a^{\eta})\ne
0$, and $\LM (a^{\gamma}a^{\beta}a^{\eta})\not\in \{ 0,1\}$, then
$\LM (a^{\gamma}a^{\alpha}a^{\eta})\prec\LM
(a^{\gamma}a^{\beta}a^{\eta})$,\par}{\parindent=0pt
then $\prec$ is called a {\it monomial ordering} on $\B$ (or a
monomial ordering on $A$).} \v5

{\bf Definition 3.2} A finitely generated $K$-algebra $A=K[a_1,\ldots ,a_n]$
is called a {\it solvable polynomial algebra} if $A$ has the PBW $K$-basis $\B =\{
a^{\alpha}=a_1^{\alpha_1}\cdots a_n^{\alpha_n}~|~\alpha
=(\alpha_1,\ldots ,\alpha_n)\in\NZ^n\}$ and a monomial ordering $\prec$ on $\B$, such that
for $\lambda_{ji}\in K^*$ and  $f_{ji}\in A$,
$$\begin{array}{l} a_ja_i=\lambda_{ji}a_ia_j+f_{ji},~1\le i<j\le n,\\
\LM (f_{ji})\prec a_ia_j~\hbox{whenever}~f_{ji}\ne 0.\end{array}$$\par

By [10], every (two-sided, respectively one-sided) ideal in a solvable polynomial algebra $A$ has a finite Gr\"obner basis with respect to a given monomial ordering; in particular, for one-sided ideals there is a noncommutative Buchberger Algorithm which has been successfully implemented in the computer algebra system \textsf{Plural}.\v5

Let $\KS =K\langle X_1,X_2,X_3\rangle$ be the free $K$-algebra generated by $X=\{X_1,X_2,X_3\}$, and let $A=\KS /\langle \G\rangle$ be a generalized down-up algebra as defined in Section 2. Writing  $A=K[a_1,a_2,a_3]$, where $a_i$ is the coset represented by $X_i$ in $A$, $1\le i\le 3$, we now start on  showing the following\v5

{\bf Theorem 3.3} If $\lambda\omega\ne 0$ and deg$f(X_1)\ge 1$ in the defining relations of $A$, then $A=K[a_1,a_2,a_3]$ is a solvable polynomial algebra in the sense of Definition 3.2. \vskip 6pt

{\bf Proof} First, it follows from Proposition 2.1 that $A$ has the PBW $K$-basis
$$\B =\{ a_2^ia_1^ja_3^{\ell}~|~(i,j,\ell)\in\NZ^3\} .$$ \par
Secondly, assigning $a_1$ the degree 1, $a_2$ and $a_3$ the degree $n=$ degree $f(a_1)$, if we take the graded lexicographic ordering $a_2\prec_{grlex}a_1\prec_{grlex}a_3$ on $\B$, then the following statements hold:\par
(a) $\prec_{grlex}$ is clearly a well-ordering on $\B$, i.e., $\prec_{grlex}$ satisfies the condition (1) of Definition 3.1.\par
(b) $\prec_{grlex}$ satisfies the condition (2) and (3) of Definition 3.1. To see this, note that  $\prec_{grlex}$ compares degree first, namely if $u,v\in\B$, then $u\prec_{grlex}v$ implies deg$u<$ deg$v$. Also note that $a_1$, $a_2$, $a_3$, $a_1a_3$, $a_2a_1$, $a_2a_3\in\B$, and $\lambda\omega\ne 0$. It follows that the generators of $A$ satisfy
$$\begin{array}{l} a_3a_1=\lambda a_1a_3-\gamma a_3~\hbox{with}~\LM (\gamma a_3)=a_3\prec_{grlex}a_1a_3~\hbox{if}~\gamma\ne 0,\\
a_1a_2=\lambda a_2a_1-\gamma a_2~\hbox{with}~\LM (\gamma a_2)=a_2\prec_{grlex}a_2a_1~\hbox{if}~\gamma\ne 0, \\
a_3a_2=\omega a_2a_3-f(a_1)~\hbox{with}~
\LM (f(a_1))\prec_{grlex}a_2a_3.\end{array}$$
Thus, for any $u=a_2^{i_1}a_1^{j_1}a_3^{\ell_1}$, $v=a_2^{i_2}a_1^{j_2}a_3^{\ell_2}\in\B$, the above properties of generators give rise to
$$\begin{array}{rcl} uv&=&(a_2^{i_1}a_1^{j_1}a_3^{\ell_1})(a_2^{i_2}a_1^{j_2}a_3^{\ell_2})\\
&=&\rho_{u,v}a_2^{i_1+i_2}a_1^{j_1+j_2}a_3^{\ell_1+\ell_2}+h
~\hbox{with}~\rho_{u,v}\in K^*, ~h\in K\hbox{-span}\B ,\\
&{~}&\hbox{deg}h<~
\hbox{deg}(a_2^{i_1+i_2}a_1^{j_1+j_2}a_3^{\ell_1+\ell_2})=(i_1+i_2+\ell_1+\ell_2)n+j_1+j_2,\\
&{~}&\hbox{thereby}~\LM(h)\prec_{grlex}a_2^{i_1+i_2}a_1^{j_1+j_2}a_3^{\ell_1+\ell_2}=\LM(uv).
\end{array}$$
This enables us to verify easily that $\prec_{grlex}$ satisfies the condition (2) and (3) of Definition 3.1. \par
Summing up, $\prec_{grlex}$ is a monomial ordering on $\B$ and therefore $A$ is a solvable polynomial algebra in the sense of Definition 3.2.\QED\v5

Concerning down-up algebras (see example (e) of Section 2), it follows from [9, 4, 5] and Theorem 3.3 that we have the following corollary.\v5

{\bf Corollary 3.4} If $\lambda\omega\ne 0$ and $f(X_1)=X_1$, then all down-up algebras $A(\alpha ,\beta ,\gamma)$ with $\alpha =\lambda +\omega$ and $\beta =-\lambda\omega$ are solvable polynomial algebras. \v5

Now, combining Proposition 2.1(i) and Theorem 3.3 with [10, 16, 17] we are able to recapture the following two basic properties of generalized down-up algebras. One may compare all proofs we presented below with those given in [5].\v5

{\bf Corollary 3.5} Let $A=\KS/\langle \G\rangle$ be a generalized down-up algebra in the sense of Section 2. The following statements hold.\par

(i) $A$ has Gelfand-Kirillov dimension three. \par

(ii) If $\lambda\omega\ne 0$, then $A$ is a Noetherian domain.\par

{\bf Proof} (i) Based on Proposition 2.1(i), this is just a special case of a more general result given in [17, P.167, Example 3], where the Ufnarovski graph [27] associated to a Gr\"obner basis is employed.\par

(ii) Since $A$ is now a solvable polynomial algebra (Theorem 3.3), this follows from the classical result  that a solvable polynomial algebra is a domain and every (left, right) ideal of $A$ has a finite Gr\"obner basis [10].\v5

\section*{4. The associated graded algebras of generalized down-up algebras}

In the literature concerning  algebras similar to $U(\textsf{sl}(2,K))$, particularly in [13, 9, 5], filtered-graded techniques have been employed for establishing some important structural properties, such as global dimension, Auslander regularity, Artin-Schelter regularity, and Koszulity, etc. By using Proposition 2.1(i) and certain related results established in [22, 16, 17], in this section we first present clearly the Gr\"obner defining relations of the associated graded algebra of a generalized down-up algebra, and then we derive (or recapture) two homological properties of generalized down-up algebras. \v5

For our purpose, we start with a little generality on the naturally filtered structure and the associated graded structure for a $K$-algebra $A=\KS /\langle \G\rangle$, where $\KS=K\langle X_1,X_2,X_3\rangle$ is the free $K$-algebra generated by $X=\{ X_1,X_2,X_3\}$, and $\langle \G\rangle$ is the (two-sided) ideal of $\KS$ generated by a subset $\G\subset\KS$. Let $\mathbb{B}$ be the standard $K$-basis of $\KS$. If  $X_1,$ $X_2,$ and $X_3$ are assigned the positive degrees $n_1$, $n_2$, and $n_3$ respectively, then $\KS$ is turned into an $\NZ$-graded algebra, i.e., $\KS=\oplus_{t\in\NZ}\KS_t$ with $\KS_t=K$-span$\{ u\in \mathbb{B}~|~\hbox{deg}u=t\}$, such that $\KS_{t_1}\KS_{t_2}\subseteq \KS_{t_1+t_2}$ for all $t_1,t_2\in\NZ$. If we further consider the grading filtration $F\KS$ of $\KS$, that is $F\KS=\{ F_q\KS\}_{q\in\NZ}$ with $F_q\KS =\oplus_{t\le q}\KS_t$, then $\KS$ is turned into an $\NZ$-filtered algebra, i.e., $\KS=\cup_{q\in\NZ}F_q\KS$ and $F_{q_1}\KS F_{q_2}\KS\subseteq F_{q_1+q_2}\KS$ for all $q_1,q_2\in\NZ$. Taking this grading filtration $F\KX$ of $\KX$ into account, the algebra $A$ has the induced filtration $FA=\{F_qA\}_{q\in\NZ}$ with $F_qA=(F_q\KX +\langle\G\rangle )/\langle\G\rangle$, such that $A=\cup_{q\in\NZ}F_qA$ and $F_{q_1}AF_{q_2}A\subseteq F_{q_1+q_2}A$ for all $q_1,q_2\in \NZ$, thereby $A$ is turned into a filtered $K$-algebra. Thus, with respect to $FA$, $A$ has its associated $\NZ$-graded $K$-algebra $G(A)=\oplus_{q\in\NZ}G(A)_q$ with $G(A)_q=F_qA/F_{q-1}A$.\par

For a nonzero $f=h_0+h_1+\cdots +h_q\in\KX$ with $h_i\in\KX_i$ and $h_q\ne 0$, we write
$$\LH (f)=h_q$$
for the highest degree homogeneous element of $f$ and call $\LH (f)$ the {\it leading homogeneous element} of $f$. \v5

{\bf Proposition 4.1} Let the algebra $A=\KX/\langle\G\rangle$ and its associated graded algebra $G(A)$ be as fixed above, and let $\prec_{grlex}$ be the graded lexicographic monomial ordering on $\KX =K\langle X_1,X_2,X_3\rangle$ with respect to some positive degrees $n_1, n_2$, and $n_3$ assigned to $X_1, X_2$, and $X_3$ respectively. Then, \par
(i) $G(A)\cong \KX/\langle\LH (\langle\G \rangle )\rangle$, where $\LH (\langle\G \rangle )=\{\LH (f)~|~f\in\langle\G\rangle\}$.\par
(ii) With respect to $\prec_{grlex}$, $\G$ is a Gr\"obner basis of $\langle\G\rangle$ if and only if  $\LH (\langle\G \rangle )$ is a homogeneous Gr\"obner basis of the graded ideal $\langle\LH (\langle\G \rangle )\rangle$. \vskip 6pt

{\bf Proof} The assertions (i) and (ii) are just special cases of [22] and [16, CH.III, CH.IV],  or more precisely, the special cases of [17, Ch.2, Theorem 3.2; Ch.4, Proposition 2.2].\QED\v5

Now, turning to generalized down-up algebras, it follows from Proposition 2.1 and Proposition 4.1 that  we have the following\v5

{\bf Theorem  4.2} With notation as in Section 2, let $A=\KX /\langle\G\rangle$ be a generalized down-up algebra and suppose that for $f(X_1)$ in $g_{32}\in\G$, the degree deg$f(X_1)\ge 1$. The following statements hold.\par
(i) If deg$f(X_1)=n\le 2$, say $f(X_1)=aX_1^2+bX_1+c$, $X_1$, $X_2$, and $X_3$ are all assigned  the degree $1$, then $\LH (\G )=\{X_3X_1-\lambda X_1X_3,~X_1X_2-\lambda X_2X_1,~
X_3X_2-\omega X_2X_3+aX_1^2\}$  is a homogeneous Gr\"obner basis of the ideal $\langle\LH (\G )\rangle$ with respect to the monomial ordering $X_2\prec_{grlex}X_1\prec_{grlex}X_3$ on $\KS$, such that $\LM (\LH (\G ))=\{ X_3X_1,~X_1X_2,~X_3X_2\}$.\par
If deg$f(X_1)=n\ge 1$,  $X_1$ is assigned the degree 1,  but $X_2$ and $X_3$ are all assigned  the degree $n$, then $\LH (\G )=\{X_3X_1-\lambda X_1X_3,~X_1X_2-\lambda X_2X_1,~
X_3X_2-\omega X_2X_3\}$  is a homogeneous Gr\"obner basis of the ideal $\langle\LH (\G )\rangle$ with respect to the monomial ordering $X_2\prec_{grlex}X_1\prec_{grlex}X_3$ on $\KS$, such that $\LM (\LH (\G ))=\{ X_3X_1,~X_1X_2,~X_3X_2\}$.\par
(ii) Let $G(A)$ be the associated graded algebra $G(A)$ of $A$ as described above, then $G(A)\cong\KX/\langle\LH(\G)\rangle$.  \par\QED\v5

{\bf Theorem 4.3} With notation and assumption as in Theorem 4.2, the following statements hold.\par
(i) $G(A)$ has homological global dimension gl.dim$G(A)=3$, $A$ has gl.dim$A\le 3$ .\par
(ii) If $\lambda\omega\ne 0$ and $f(X_1)=aX_1^2+bX_1+c$ with $ab\ne 0$, then $G(A)$ is the associated graded algebra of the conformal $\textsf{sl}(2,K)$ enveloping algebra, thereby all results obtained in [13] hold true for $G(A)$, particularly $G(A)$ is an Auslander regular algebra satisfying Cohen-Macaulay property (see the definitions given in Section 1), and so too is $A$.\par

(iii) If $\lambda\omega\ne 0$ and deg$f(X_1)=n\ge 1$, then $G(A)$ is a solvable polynomial algebra in the sense of Definition 3.2. \par

(iv) If deg$f(X_1)=n\ge 1$, then, with $\lambda\omega\ne 0$, $G(A)$ is an Auslander regular algebra satisfying Cohen-Macaulay property, and so too is $A$.\vskip 6pt

{\bf Proof} (i) This follows from Theorem 4.2 and [17, Ch.5, Corollary 7.6].\par
(ii) By Theorem 4.2, if deg$f(X_1)=n\le 2$, $X_1$, $X_2$, and $X_3$ are all assigned  the degree $1$, then $G(A)\cong \KX/\langle\LH(\G)\rangle$
where $$\LH(\G)=\{X_3X_1-\lambda X_1X_3,~X_1X_2-\lambda X_2X_1,~
X_3X_2-\omega X_2X_3+aX_1^2\}.$$ Note that if $\lambda\omega a\ne 0$, then $G(A)$ is  a conformal $\textsf{sl}(2,K)$ enveloping algebra, thereby all results obtained in [13] hold true for $A$, particularly $G(A)$ is an Auslander regular algebra satisfying Cohen-Macaulay property, and so too is $A$.\par

(iii) By Theorem 3.3, $A=K[a_1,a_2,a_3]$ is a solvable polynomial algebra with respect to the graded lexicographic monomial ordering $a_2\prec_{grlex}a_1\prec_{grlex}a_3$, where $a_i$ is the coset represented by $X_i$ in $A$ ($=\KS /\langle\G \rangle$), and $a_1$ is assigned the degree 1 and $a_2$, $a_3$ are all assigned the degree $n=$ deg$f(X_1)$. So the proof of this assertion is just an analogue of the proof of a similar result  established in [21, Section 3] and [16, CH.IV, Section 4] (where general quadric solvable polynomial algebras are considered).\par

(iv) By Theorem 4.2, if deg$f(X_1)=n\ge 1$, $X_1$ is assigned the degree 1, but $X_2$ and  $X_3$ are all assigned  the degree $n$, then $G(A)\cong \KX/\langle\LH(\G)\rangle$
where $\LH(\G)=\{X_3X_1-\lambda X_1X_3,~X_1X_2-\lambda X_2X_1,~
X_3X_2-\omega X_2X_3\}$ with $\LM (\LH (\G))=\{ X_3X_1, X_1X_2, X_3X_2\}$. Note that if $\lambda\omega\ne 0$, then $G(A)$ is clearly a skew polynomial algebra over $K$. It follows from [20] (or by a similar argumentation as given in [13, 18, 6]) that $G(A)$ is an Auslander regular algebra satisfying Cohen-Macaulay property, and so too is $A$.\v5

\section*{5. The homogenized algebras of generalized down-up algebras}
Let $\KX =K\langle X_1,X_2,X_3\rangle$ be the free $K$-algebra generated by $X=\{ X_1,X_2,X_3\}$ and consider a generalized down-up algebra $A=\KX/\langle\G\rangle$ in the sense of Section 2.
In this section, by using the Gr\"obner basis techniques we establish some basic structural properties of the homogenized algebras of generalized down-up algebras (see the definition  below).\v5

We start by  recalling the general  notion of a homogenized algebra and some relevant results from the literature (for instance, [25, 13, 15, 14, 22, 16, 17]).\v5

Let the free $K$-algebra $\KX =K\langle X_1,\ldots ,X_n\rangle$ be equipped with the graded structure $\KX=\oplus_{q\in\NZ}\KX_q$ by assigning each generator $X_i$ a positive degree  deg$X_i=m_i>0$, $1\le i \le n$. Considering the free $K$-algebra $K\langle X,T\rangle =K\langle X_1,\ldots ,X_n,T\rangle$ generated by $\{ X_1,\ldots ,X_n,T\}$, let $K\langle X,T\rangle$ be equipped with the the graded structure $K\langle X,T\rangle =\oplus_{q\in\NZ}K\langle X,T\rangle_q$ by assigning deg$X_i=m_i$ for $1\le i \le n$, and deg$T=1$. Then,{\parindent=1.4truecm\par

\item{(a)} for a nonzero element $f=\sum_{i=1}^mh_i\in\KX$ with $h_i\in\KX_{q_i}$, $q_1<q_2<\ldots <q_m$, and $h_m\ne 0$, the degree-$q_m$ homogeneous element $\T{f}=\sum_{i=1}^mT^{q_m-q_i}h_i$ of $K\langle X,T\rangle$ is referred to as the {\it homogenization of $f$} in $K\langle X,T\rangle$ with respect to $T$;\index{homogenization}\par

\item{(b)} for a nonempty subset $S\subset \KX$, writing $I=\langle S\rangle$ for the ideal generated by $S$, $A=\KX /I$, and writing $\T{S}=\{\T{f}~|~f\in S\}\cup\{X_iT-TX_i~|~1\le i\le n\}$, the quotient algebra $H(A)=K\langle X,T\rangle /\langle\T{S}\rangle$ is referred to as the {\it homogenized algebra} of $A$  with respect to $T$.\par}{\parindent=0pt\v5

{\bf Remark} Related to the theory of quantum groups, the study of homogenized enveloping algebra was  proposed by S.P. Smith in [25], and the study of homogenized down-up algebra was proposed by G. Benkart and T. Roby in [2].  In [13], homogenized conformal $\textsf{sl}_2$ enveloping algebras were used to study modules over conformal $\textsf{sl}_2$ enveloping algebras. Concerning the representation theory and noncommutative geometry of homogenized enveloping algebras, the reader is referred to, for instance, [14, 15, 13]. A study of homogenized algebras, via  Gr\"obner  defining relations of algebras, was initiated in [22, part 2] and further extended in [16, 17]. }\v5

Furthermore, consider the $\NZ$-grading filtration $F\KX =\{ F_q\KX\}_{q\in\NZ}$ of $\KX$ with each $F_q\KX =\oplus_{\ell\le q}\KX_{\ell}$, which gives rise to the $\NZ$-filtration $FA =\{ F_qA\}_{q\in\NZ}$ of the algebra $A=\KX /I$ with each $F_qA=F_q\KX +I/I$. Since $F_{q_1}AF_{q_2}A\subseteq F_{q_1+q_2}A$ for all $q_1,q_2\in\NZ$, this naturally gives rise to an $\NZ$-graded algebra $\T A=\oplus_{q\in\NZ}\T{A}_q$ with $\T{A}_q =F_qA$. This graded algebra is usually referred to as the {\it Rees algebra} of $A$ determined by $FA$ (see [20], or [16, 17] for more detailed discussion on $\T{A}$). Concerning the algebras $H(A)$ and $\T{A}$, the following proposition stem from [22] and [17, Section 7.2].
\v5

{\bf Proposition 5.1} With notation and the degrees assigned to the $X_i$ and $T$
as fixed above, let $\langle\T{I}\rangle$ be the ideal of $K\langle X,T\rangle$ generated by $\T{I}=\{\T{f}~|~f\in I\}\cup\{X_iT-TX_i~|~1\le i\le n\}$. Taking  a graded lexicographic ordering $\prec_{grlex}$  on the standard $K$-basis $\mathbb{B}$ of $\KX$ such that
$$X_{i_1}\prec_{grlex}X_{i_2}\prec_{grlex}\cdots\prec_{grlex}X_{i_n},\quad 1\le i\le n,$$
and extending $\prec_{grlex}$ to the graded lexicographic ordering $\prec_{_{T\hbox{\tiny -}grlex}}$ on the standard $K$-basis $\mathbb{B} (T)$ of $K\langle X,T\rangle$ such that
$$T\prec_{_{T\hbox{\tiny -}grlex}}X_{i_1}\prec_{_{T\hbox{\tiny -}grlex}}X_{i_2}\prec_{_{T\hbox{\tiny -}grlex}}\cdots\prec_{_{T\hbox{\tiny -}grlex}}X_{i_n},\quad 1\le i\le n,$$
if $\G$ is a Gr\"obner basis of $I=\langle S\rangle$ with respect to $\prec_{grlex}$, then
$$\T{\G}=\{\T{g}~|~g\in \G\}\cup\{X_iT-TX_i~|~1\le i\le n\}$$
is a homogeneous Gr\"obner basis of the graded ideal $\langle\T{I}\rangle$ in $K\langle X,T\rangle$ with respect to $\prec_{_{T\hbox{\tiny -}grlex}}$, $\langle\T{S}\rangle =\langle\T{I}\rangle$, and there is a graded algebra isomorphism
$$H(A)=K\langle X,T\rangle /\langle\T{\G}\rangle\cong\T{A}.$$\par\QED\v5

Applying Proposition 5.1 to generalized down-up algebras, we may derive the following\v5

{\bf Theorem 5.2} Let $A=\KX/\langle\G\rangle$ be a generalized down-up algebra with def$f(X_1)=n\ge 1$ in $\G$. With notation as in Section 2 and those fixed above, the following statements hold.\par
(i) $\T{\G}=\{\T{g}~|~g\in\G\}\cup\{ X_iT-TX_i~|~1\le i\le 3\}$ is a homogeneous Gr\"obner basis of the ideal $\langle\T{\G}\rangle$ in the free $K$-algebra $K\langle X, T\rangle =K\langle X_1,X_2,X_3,T\rangle$ with respect to the graded lexicographic monomial ordering $$T\prec_{_{T\hbox{\tiny -}grlex}}X_2\prec_{_{T\hbox{\tiny -}grlex}}X_1\prec_{_{T\hbox{\tiny -}grlex}}X_3$$ on $K\langle X,T\rangle$, where
$$\hbox{deg}T=1=~\hbox{deg}X_1,\quad \hbox{deg}X_2=~\hbox{deg}X_3=n=~\hbox{deg}f(X_1).$$\par

(ii) The homogenized algebra $H(A)$ of $A$ is isomorphic to the Rees algebra of $A$, i.e.,  $H(A)=K\langle X,T\rangle /\langle\T{\G}\rangle\cong \T{A}$, where the filtration $FA$ of $A$ is the one induced by the grading filtration of $K\langle X,T\rangle$ (degrees of $T, X_1,X_2,$ and $X_3$ are assigned as in (i)).\vskip 6pt

{\bf Proof} (i) By Proposition 2.1, the set $\G$ of defining relations of $A$ forms a Gr\"obner basis of the ideal $\langle\G\rangle$ with respect to the graded lexicographic monomial ordering $$X_2\prec_{grlex}X_1\prec_{grlex}X_3,$$
where deg$X_1=1$, deg$X_2=$ deg$X_3=n=$ deg$f(X_1)\ge 1$. It follows from Proposition 5.1 that the assertions (i) and (ii) hold true.\QED\v5

{\bf Theorem 5.3} Let $A=\KX/\langle\G\rangle$ be a generalized down-up algebra with deg$f(X_1)=n\ge 1$ in $\G$. With notation as made above, the following statements hold.\par

(i) The Gelfand-Kirillov dimension of $\T{A}$ and $H(A)$ are all equal to 4.\par

(ii) gldim$\T{A}=$ gl.dim$H(A)=4$.\par

(iii) The Hilbert series of $\T{A}$ and $H(A)$ is $\frac{1}{(1-t)^4}$. {\parindent=0pt\par
In particular, the assertions (i), (ii), and (iii) hold true for all down-up algebras.}\par

(iv) If $f(X_1)=aX_1^2+bX_1+c$, then $\T{A}$ is a classical quadratic Koszul algebra, and so too is $H(A)$. \par

Moreover, if $\lambda\omega\ne 0$ in $\G$, then the following statements hold.\par

(v) $\T{A}$ is an $\NZ$-graded solvable polynomial algebra in the sense of Definition 3.2, and so too is $H(A)$.\par

(vi) $\T{A}$ is a Noetherian domain, and so too is $H(A)$.\par

(vii) $\T{A}$ is an Auslander regular algebra satisfying the Cohen-Macaulay property, and so too is $H(A)$.
\vskip 6pt

{\bf Proof} Due to the graded algebra isomorphism $\T{A}\cong H(A)$ (Proposition 5.1), it is clear that we need only to prove all assertions for $\T{A}$ below.\par
By Theorem 5.2, taking the monomial ordering $$T\prec_{_{T\hbox{\tiny -}grlex}}X_2\prec_{_{T\hbox{\tiny -}grlex}}X_1\prec_{_{T\hbox{\tiny -}grlex}}X_3$$ on $K\langle X,T\rangle$, where $\hbox{deg}T=1=~\hbox{deg}X_1,\quad \hbox{deg}X_2=~\hbox{deg}X_3=n=~\hbox{deg}f(X_1),$ the Gr\"obner basis $\T{\G}$ of $\langle\T{\G}\rangle$ is now consisting of
$$\begin{array}{l} \T{g}_{31}=X_3X_1-\lambda X_1X_3+\gamma TX_3,\\
\T{g}_{12}=X_1X_2-\lambda X_2X_1+\gamma TX_3,\\
\T{g}_{32}=X_3X_2-\omega X_2X_3+\sum^n_{i=0}a_iT^{2n-i}X_1^i~(\hbox{provided}~f(X_1)=\sum^n_{i=0}a_iX_1^i),\\
X_iT-TX_i,~1\le i\le 3,\end{array}$$
such that $\LM (\T{\G})=\{ X_3X_1,~X_1X_2,~X_3X_2, ~X_1T,~X_2T,~X_3T\}$. So, the assertions (i), (ii), and (iii) follow from [17, Ch.7, Corollary 3.6].\par

(iv) If $f(X_1)=aX_1^2+bX_1+c$, then by Proposition 2.1 (or the argumentation made before it) we know that in this case the Gr\"obner basis  $\G$ is obtained with respect to the monomial ordering $X_2\prec_{grlex}X_1\prec_{grlex}X_3$, but $X_1, X_2, X_3$ are all assigned the natural degree 1. Thus, $\T{\G}$ consists of quadratic relations. Hence $\T{A}$ is a classical quadratic Koszul algebra. \par

We next proceed to deal with the case that $\lambda\omega\ne 0$ in the Gr\"obner basis $\G$ of $A$.\par

(v) By Theorem 3.3, $A=K[a_1,a_2,a_3]$ is a solvable polynomial algebra with respect to the graded lexicographic monomial ordering $a_2\prec_{grlex}a_1\prec_{grlex}a_3$, where $a_i$ is the coset represented by $X_i$ in $A$,  $a_1$ is assigned the degree 1, and $a_2$, $a_3$ are all assigned the degree $n=$ deg$f(X_1)$. So the proof of this assertion is just an analogue of the proof of a  result  established in [21, Section 3] and [16, CH.IV, Section 4] (where  general quadric solvable polynomial algebras are considered).\par

(vi) As in the proof of Corollary 3.5, since $\T{A}$ is now a solvable polynomial algebra by the assertion of (v), this follows from the classical result that a solvable polynomial algebra is a domain and every (left, right) ideal of A has a finite Gr\"obner basis [10].\par

(vii) By Theorem 4.3(iii), $A$ and its associated graded algebra $G(A)$ are Auslander regular algebra satisfying Cohen-Macaulay property. It follows from [19], [20, Ch.III], and [16, Ch.III]) that $\T{A}$ is an Auslander regular algebra satisfying Cohen-Macaulay property, and so too is $H(A)$. \QED\v5

{\bf Remark} In [5], it was already shown, in a quite different way, that if  $\lambda\omega\ne0$, then $H(A)$ is a Noetherian domain and an Artin Achelter regular algebra of global dimension 4.\v5

\centerline{Refeerence}{\parindent=1.47truecm\par

\item{[1]}  Benkart, G. (1999). Down-up algebras and Witten's deformations of the
universal enveloping algebra of sl2. {\it Contemp. Math}. 224:29--45.

\item{[2]} Benkart, G.,  Roby, T. (1999). Down-up algebras. {\it J. Alg}. 209:305--
344. Addendum: {\it J. Alg}. 213: 378.

\item{[3]} Carvalho, Paula A. A. B.,   Lopes, Samuel A. (2009). Automorphisms of Generalized Down-Up Algebras. {\it Comm. Algebra}. 37(5):1622--1646.

\item{[4]}  Carvalho, P.  Musson, M. (2000). Down-up algebras and their representation
theory. {\it J. Alg}. 228:286--310.

\item{[5]} Cassidy, T., Shelton, B. (2004). Basic properties of generalized down-up
algebras. {\it J. Alg}. 279:402--421.

\item{[6]} Cassidy, T., Shelton, B. (2007). PBW-deformation theory and regular
central extensions. {\it J. f\"ur die reine und angewandte Mathematik}.
610:1--12.

\item{[7]}  Dubois-Violette, M.,  Popov, T. (2002). Homogeneous algebras, statistics and combinatorics. {\it Letters in Mathematical Physics}. 61:159--170.

\item{[8]} Jie, X., Van Oystaeyen, F. (1995). Weight modules and their extensions over a class of algebras similar to the enveloping algebra of sl(2,C). {\it J. Algebra}.
    175:844--864.

\item{[9]}  Kirkman, E.,  Musson, I. M.,  Passman, D. S. (1999). Noetherian down-up
algebras. {\it Proc. Amer. Math. Soc}. 127:2821--2827.

\item{[10]}  Kandri-Rody, A., Weispfenning, V. (1990).  Non-commutative Gr\"obner
bases in algebras of solvable type. {\it J. Symbolic Comput}. 9:1--26.

\item{[11]}  Kulkarni, Rajesh S. (2001). Down-Up algebras and their representations. {\it J. Algebra}. 245:431--462.

\item{[12]} Le Bruyn, L. (1994). two remarks on Witten's quantum sl2 enveloping algebras.
{\it Comm. Alg}. 22:865--876.

\item{[13]}  Le Bruyn, L. (1995). Conformal sl2 enveloping algebras. {\it Comm. Alg}.  23:1325--1362.

\item{[14]}  Le Bruyn, L.,  Smith, S. P. (1993). Homogenized
sl(2). {\it Proc. Amer. Math. Soc}. 3(118):725--730.

\item{[15]} Le Bruyn, L., Van den Bergh, M. (1993). On quantum spaces of Lie
algebras. {\it Proc. Amer. Math. Soc.} 2(119):407--414.

\item{[16]} Li, H. (2002). {\it Noncommutative Gr\"obner Bases and Filtered-graded Transfer}. Lecture Notes in Mathematics, Vol. 1795. Springer, Berlin, Heidelberg. \url{https://doi.org/10.1007/978-3-540-45765-7_10 }

\item{[17]} Li, H. (2011). {\it Gr\"obner Bases in Ring Theory}. World Scientific Publishing Co. \url{https://doi.org/10.1142/8223 }

\item{[18]}  Le  Bruyn, L.,  Smith, S. P.,    Van  den  Bergh, M. (1996).  Central  extensions  of three  dimensional  Artin-Schelter regular  algebras. {\it Math. Z}.  222:171-212.

\item{[19]} Li, H., Van Oystaeyen, F. (1991). Global dimension and Auslander regularity
of Rees rings. {\it Bull. Soc. Math. Belg.} 43:59--87.

\item{[20]} Li, H., Van Oystaeyen, F. (1996, 2003). {\it Zariskian Filtrations}. K-Monograph
in Mathematics, Vol.2. Kluwer Academic Publishers, Berlin Heidelberg: Springer-Verlag.

\item{[21]} Li, H., Wu,Y. (2000). Filtered-graded transfer of
Gr\"obner basis computation in solvable polynomial algebras. {\it
Communications in Algebra}. 1(28):15--32.

\item{[22]} Li, H., Wu, Y., Zhang, J. (1999). Two applications of
noncommutative Gr\"obner bases. {\it Z. Ann. Univ. Ferrara}. 1(45):1--24.

\item{[23]} Mora, T. (1994).  An introduction to commutative and noncommutative
Gr\"obner Bases. {\it Theoretic Computer Science}. 134:131--173.

\item{[24]} Smith, S. P. (1990). A class of algebras similar to the enveloping algebra of $sl(2)$. {\it Trans. Amer. Math. Soc}. 332:285--314.

\item{[25]} Smith, S. P. (1992). Quantum groups: an introduction and survey
for ring theorists. In: Montgomery, S., Small, L., ed(s). {\it Noncommutative Rings},  MSRI Publ.  New York, Springer-Verlag, pp. 131--178.

\item{[26]} Woronowicz, S. L. (1987). Twisted SU(2) group. An example of a noncommutative
differential calculus. {\it Publ. RIMS} 23:117--181.

\item{[27]} Ufnarovski, V. (1982). A growth criterion for graphs and algebras defined
by words. {\it Mat. Zametki}, 31:465--472 (in Russian); English
translation: {\it Math. Notes}, 37:238--241.

\item{[28]} Witten, E. (1990). Gauge theories, vertex models, and quantum groups. {\it Nuclear Phys}. B 330:285--346.
\end{document}